\documentclass{amsart}

\usepackage{amssymb}

\usepackage[all]{xy}

\newtheorem{thm}{Theorem}%[section]

\newtheorem{lem}[thm]{Lemma}

\newtheorem{prop}[thm]{Proposition}

\theoremstyle{definition}
\newtheorem{defn}[thm]{Definition}

\newtheorem{say}[thm]{}
\newtheorem{exmp}[thm]{Example}

\newtheorem{ques}[thm]{Question}    %!!!!!!!!!!!!!!!!!!!!

\newtheorem*{ack}{Acknowledgments}      % \renewcommand{\theack}{} 

\newtheorem{defn-thm}[thm]{Definition--Theorem}  %!!!!!!!!!!!!!!!!!!!!!!!!
\newtheorem{defn-lem}[thm]{Definition--Lemma}  %!!!!!!!!!!!!!!!!!!!!!!!!
\theoremstyle{remark}

\renewcommand{\c}[0]{{\mathbb C}}  

\renewcommand{\o}[0]{{\mathcal O}} 

\newcommand{\n}[0]{{\mathbb N}}
\newcommand{\p}[0]{{\mathbb P}}

\newcommand{\qtq}[1]{\quad\mbox{#1}\quad}

\newcommand{\pic}[0]{\operatorname{Pic}}

\newcommand{\rank}[0]{\operatorname{rank}}

\newcommand{\codim}[0]{\operatorname{codim}}    
\newcommand{\im}[0]{\operatorname{im}}

\newcommand{\sing}[0]{\operatorname{Sing}}

\newcommand{\nec}[1]{\overline{NE}({#1})}

\newcommand{\tsum}[0]{\textstyle{\sum}}

\newcommand{\GL}{\mathrm{GL}}

\newcommand{\SU}{\mathrm{SU}}

\def\into{\DOTSB\lhook\joinrel\to}

\def\loccoh#1.#2.#3.#4.{H^{#1}_{#2}(#3,#4)}

\DeclareMathAlphabet{\mathchanc}{OT1}{pzc}%
                                {m}{it}

\newcommand{\chain}[0]{\operatorname{Chain}}

\begin{document}
\bibliographystyle{amsalpha}

%\today

\title[Neighborhoods in  homogeneous spaces]{Neighborhoods of  subvarieties 
in  \\ homogeneous spaces}
\author{J\'anos Koll\'ar}

\maketitle

Let $X$ be an algebraic  variety over $\c$
and $D\subset X$ a Euclidean open subset. 
It is interesting to find connections between the function theory or topology
of $D$ and $X$.
There is not much to say if $D$ is affine or Stein.
By contrast, strong results are known 
if $D$ contains a positive dimensional, compact subvariety $Z$
with ample normal bundle:
\begin{itemize}
\item The field of meromorphic functions ${\mathcal Mer}(D)$ 
is a finite extension of the field of rational functions
${\mathcal Rat}(X)$. The proof, by \cite{MR0241433, MR0232780},
 relies on cohomology vanishing for symmetric powers of the
normal bundle of $Z$. 
\item  The image of the natural map
$\pi_1(D)\to \pi_1(X)$ has finite index in $\pi_1(X)$.
More generally, for every Zariski open subset $X^0\subset X$,
the image of the map
$\pi_1(D\cap X^0)\to \pi_1(X^0)$ has finite index in $\pi_1(X^0)$.
The proof, by \cite{nap-ram}, uses $L^2$ $\bar\partial$-methods.
\end{itemize}

%${\mathcal Mer}(D)$ on $D$ and the field of rational functions
%${\mathcal Rat}(X)$ on $X$.

The isomorphism of these function fields
and the surjectivity of the maps between the fundamental groups are 
subtler questions. 
 ${\mathcal Mer}(D)={\mathcal Rat}(X)$ was proved  for
$\p^n$   \cite{MR0241433,  MR0251043} and for
Grassmannians  \cite{MR650388}. 
The surjectivity of the maps between the fundamental groups was
established for neighborhoods of certain high degree rational curves
in \cite{MR1786496, MR2011744}.

It was also observed by \cite{MR0241433}
that if ${\mathcal Mer}(D)={\mathcal Rat}(X)$
for every  $Z$ and $D$  then $X$ is 
simply connected, but the close connection between the two types of
theorems
was not fully appreciated.

I was lead to consider these topics while trying to answer 
some  problems about non-classical flag domains   raised by
Griffiths and  Toledo during the conference
{\it  Hodge Theory and Classical Algebraic Geometry;}
see Question \ref{grt.ques}.

It turns out that the answer needs very few properties of
non-classical flag domains. The natural setting is
to study an arbitrary, simply connected,
quasi projective, homogeneous space $X$,
a proper subvariety $Z\subset X$ and a Euclidean open neighborhood
$D\supset Z$. Theorem \ref{main.thm} gives a complete
description of those pairs $Z\subset X$ for which the
holomorphic/meromorphic function theory of $D$
is determined by the regular/rational function theory of $X$.
The precise connection is established through an 
 understanding of the surjectivity of $\pi_1(D\cap X^0)\to \pi_1(X^0)$.

We allow $Z$ to be
singular and with non-ample normal sheaf.
A slight difference is that, while \cite{MR0241433, MR0232780,  MR0251043}
studied the formal completion of $X$ along $Z$, we work with
actual open neighborhoods. In the ample normal bundle case
the two versions are equivalent, but I am not sure
that this also holds in general; cf.\ \cite{MR0206980}.

The main tool is the study of  chains made up of
translates of $Z$  in $X$ and in $D$. 
In the projective case such techniques form  the basis of the
study of rationally connected varieties; see \cite{rc-book}
for  a detailed treatment  or \cite{ar-ko} for more  introductory
lectures. For non-proper homogeneous spaces these ideas were
used in \cite{MR1369412}.

\begin{defn}
\label{nondeg.defn}
Let $X=G/H$ be a  simply connected, quasi projective, homogeneous space.
The left action of $g\in G$ on $X$ is denoted by  $\tau_g$; we call it a
{\it translation.}

An irreducible subvariety $Z\subset X$
 will be called 
{\it degenerate} if there is a subgroup
$H\subset K\subset G$ such that  $Z$ is contained in a fiber of
the natural projection $p_K:G/H\to G/K$;  otherwise we 
call $Z$ {\it nondegenerate.}
(If $X$ is not simply connected,  these  notions 
should be modified; see  Example \ref{H.disconn.exmp}.)

For example, if $X$ is a projective homogeneous space
of Picard number 1 then every positive dimensional
subvariety is nondegenerate.
More generally, if the $X_i$ are projective homogeneous spaces
of Picard number 1 then $Z\subset \prod X_i$ is 
nondegenerate iff none of the coordinate projections
$Z\to X_i$ is constant. 
\end{defn}

Our main theorem is the following.

\begin{thm}\label{main.thm} Let $G$ be a connected algebraic group over $\c$,
 $X=G/H$  a quasi projective, simply connected, homogeneous space
 and $Z\subset X$  a compact, irreducible subvariety.
Let $D\subset X$ denote a sufficiently small 
Euclidean  open neighborhood of $Z$.
The following are equivalent.
\smallskip

\noindent (Finiteness conditions)
\begin{enumerate}
\item $H^0(D, \o_D)=\c$.
\item  $H^0(D, L)$ is finite dimensional for every line bundle $L$ on $D$.
\item  $\dim H^0(D, L^m)=O\bigl(m^{\dim D}\bigr)$ 
 for every line bundle $L$ on $D$.
\item  $H^0(D, E)$ is finite dimensional for every 
coherent, torsion free sheaf $E$.
\end{enumerate}
(Isomorphism conditions)
\begin{enumerate}\setcounter{enumi}{4}
\item  $H^0(D, L|_D)\cong H^0(X, L)$  for every line bundle $L$ on $X$.
\item  $ H^0(D, F|_D)\cong H^0(X, F)$  for every coherent, reflexive 
sheaf $F$ on $X$.
\item ${\mathcal Mer}(D)={\mathcal Rat}(X)$.
\item The conditions {\rm (1--7)} hold for every
finite, \'etale cover $\tilde D\to D$.
\end{enumerate}
(Fundamental group conditions on Zariski open subsets  $X^0\subset X$)
\begin{enumerate}\setcounter{enumi}{8}
\item $\pi_1\bigl(D\cap X^0\bigr)\to \pi_1\bigl(X^0\bigr)$ 
is surjective for every $X^0$.
\item $\pi_1\bigl(u_D^{-1}(X^0)\bigr)\to \pi_1\bigl(X^0\bigr)$ 
is surjective  for every $X^0$ and
 every finite, \'etale cover $u_D:\tilde D\to D$.
\item $\pi_1\bigl(\tau_g(Z)\cap X^0\bigr)\to \pi_1\bigl(X^0\bigr)$
is surjective   for every $X^0$ and  general $g\in G$.
\item $\pi_1\bigl((\tau_g\circ u_Z)^{-1} (X^0)\bigr)\to \pi_1\bigl(X^0\bigr) $
is surjective  for every $X^0$,
 every finite cover $u_Z:\tilde Z\to Z$ and  general $g\in G$.
\end{enumerate}
(Geometric characterizations of $Z$)
\begin{enumerate}\setcounter{enumi}{12}
\item $Z\cap B\neq \emptyset$ for every nonzero divisor $B\subset X$.
\item For every $x_1, x_2\in X$ there is a connected subvariety
$Z(x_1, x_2)\subset X$ containing them, whose irreducible components
are  translates of $Z$.
\item Same as {\rm (14)} with at most $ 2\dim X$  irreducible components.
\item $Z$ is nondegenerate in $X$.
\end{enumerate}
\end{thm}

\begin{say}[Comments]  We will show that 
(\ref{main.thm}.13--16) $\Rightarrow$ (\ref{main.thm}.1--12) 
for every $D$. The precise conditions for the other
implications vary.  
In all cases of (\ref{main.thm}.1--7),  the space of global sections
gets bigger as $D$ gets smaller.
For (\ref{main.thm}.9--10) 
the relevant assumption is that $D$ retracts to $Z$, or  that it is
contained in a neighborhood that retracts to $Z$.

Many parts of Theorem \ref{main.thm} work even if $X$ is not
simply connected, but the deepest statements, 
(\ref{main.thm}.5--12)  do not.
In one of the most interesting cases, when $Z$ is a smooth, rational curve,
there are simply connected neighborhoods $D\supset Z$. 
Thus $D\into X$ lifts to the universal cover  $D\into \tilde X$,
hence the function theory of $D$ is determined by $\tilde X$;
the embedding  $D\into X$ is just an accident.

The finite dimensionality statements
(\ref{main.thm}.1--4)  fit in the general framework of the papers
\cite{MR0241433, MR0232780,  MR0251043}.

 The isomorphism statements
(\ref{main.thm}.5--7) are more subtle.  They were known for
$\p^n$   \cite{MR0241433,  MR0251043} and for
Grassmannians  \cite{MR650388}.
In the terminology of \cite{MR0241433}, property (\ref{main.thm}.7)
is called the  G3 condition. It has been investigated in many other cases,
see  \cite{MR0332788, MR556311, MR2555947, MR2876923}.

Condition (\ref{main.thm}.8)  
mixes together some obvious claims with some 
quite counter intuitive ones.
If $v:D'\to D$ is a finite (possibly ramified) cover
and  $E'$ is a coherent, torsion free sheaf on $D'$
then $v_*E'$ is also  a coherent, torsion free sheaf
and $H^0(D', E')= H^0(D, v_*E')$. 
Thus (\ref{main.thm}.2--4)  hold for $D'$ as well. 
By contrast, one would expect to find more
sections and meromorphic functions on $\tilde D$.
In particular,  (\ref{main.thm}.8)   implies that a nontrivial 
finite \'etale cover $\tilde D\to D$  is never embeddable
into any algebraic variety.
%see Paragraph \ref{covers.of.D.say} for more details.

A weaker version of the Lefschetz--type properties (\ref{main.thm}.9--12),
asserting finite index image instead of surjectivity,
is roughly equivalent to the finiteness (\ref{main.thm}.2); see
 \cite{nap-ram}. The stronger variants are
studied in the papers  \cite{MR1786496, MR2011744}
 when $Z$ is a  rational curve. 
In (\ref{main.thm}.11--12) the adjective {\it general}
means that the claim holds for all $g$ in a nonempty Zariski open subset  
$U(X^0)\subset G$ which depends on $X^0$.
Earlier results 
gave  (\ref{main.thm}.11)  for sufficiently
high degree curves only.

The stronger forms (\ref{main.thm}.10) and  (\ref{main.thm}.12)
 may seem surprising at first
since by taking \'etale covers, the groups
$\pi_1\bigl( u_D^{-1} (X^0)\bigr)$ are getting smaller.
However, $X$ itself is simply connected, thus
all the fundamental group of $X^0$ comes from
loops around  $X\setminus X^0$, and  such loops 
are preserved by \'etale covers of $D$.

Presumably (\ref{main.thm}.15) also holds with 
at most $ \dim X$  irreducible components
or maybe with some even smaller linear function of $\dim X$.

If $X$ is projective, then the various projections
$G/K\to G/H$ correspond to the faces of the cone of curves
$\nec{X}$. Thus  a curve
$C\subset X$ is nondegenerate iff its homology class
$[C]\in \nec{X}$ is an interior point.

%\end{say} 
%\begin{say}[Generalizations] 

Many of the conditions in 
Theorem \ref{main.thm} are equivalent to each other under
much more general conditions.

One key assumption for an arbitrary pair
$Z\subset X$ is that the deformation theory of $Z$ in $X$ should be
as rich as for homogeneous spaces. Under such conditions,
the properties within any of the 4 groups tend to be equivalent to
each other.  However, I could not write down neat, general
versions in all cases.

A rather subtle point is the role of the simple connectedness of $X$.
While this is definitely needed, it seems  more important to
know that the stabilizer subgroup $H$ is connected.

The equivalence of the 4 groups to each other is more complicated
and it depends on further properties. Even
 when $Z$ is a smooth rational curve with ample normal bundle,
the conditions (\ref{main.thm}.5--12) are much stronger
than (\ref{main.thm}.1--4).
% see Example \ref{leff.exmp.3}.

The latter case has been studied in the papers 
\cite{MR1369412, MR1786496, MR2011744, MR2019976}
and most of the arguments of this note have their origins
in one of them.
\end{say}

\begin{say}[Proving that (\ref{main.thm}.1--15)  
all imply (\ref{main.thm}.16)]{\ }

Assume that  $Z\subset X$ is degenerate.
Thus there is a subgroup
$H\subset K\subset G$ such that
$Z$ is contained in one of the fibers of
$p_K:G/H\to G/K$. 

 Let $U\subset G/K$ be a small Stein neighborhood of the point 
$p_K(C)$. Then $p_K^{-1}(U)\subset X$ is an 
open neighborhood of $C$ with many holomorphic functions
and (\ref{main.thm}.1--8) all fail for every neighborhood contained in 
$p_K^{-1}(U)$.

Similarly, if $U$ is contractible and 
$U\subset Y^0\subset G/K$ is Zariski open such that
$\pi_1(Y^0)$ is infinite then (\ref{main.thm}.9--12) fail for
 every Zariski open subset  of 
$p_K^{-1}(Y^0)$.

The preimage of a divisor $B_K\subset G/K$ shows that (\ref{main.thm}.13) fails
and translates of $Z$ never connect points in different fibers of $p_K$.\qed 
\end{say}

% It remains to prove that if  $Z$ is nondegenerate 
% then all the other properties hold. 

\subsection*{Open problems}{\ }
\medskip

In connection with Theorem \ref{main.thm} an interesting
open problem is to understand which (non-proper) homogeneous spaces
$X=G/H$ contain a proper, nondegenerate subvariety.
Consider the following conditions.
\begin{itemize}
\item There is a projective compactification $\bar X\supset X$
such that $\bar X\setminus X$ has codimension $\geq 2$.
\item $X$ contains a proper, nondegenerate subvariety.
\item There is no subgroup $H\subset K\subsetneq G$
such that $G/K$ is quasi affine.
\end{itemize}
It is clear that each one implies the next.
Based on  \cite[Sec.6]{MR1369412}, one can ask the following.

\begin{ques} Are the above 3 conditions equivalent for a
homogeneous space?
\end{ques} 

\begin{ques} If $X$ contains a proper, nondegenerate subvariety,
does it contain a proper, nondegenerate,
smooth, rational curve?
\end{ques} 

\subsection*{Applications to non-classical flag domains}{\ }
\medskip

Our results can be used to study global sections of coherent sheaves 
over certain homogeneous complex manifolds.
While traditionally most attention was devoted to 
compact homogeneous spaces and to 
 Hermitian symmetric domains,
other examples have also been studied
\cite{MR0251246, MR2188135}.
The recent paper \cite{grt} studies  the geometry
of {\it non-classical flag domains.} 
Most period domains of Hodge structures are of this type.
For our purposes the precise definition is
not important, we  need only two of their properties.
\begin{itemize}
\item[$\circ$] A flag domain  is an  open subset  of
a projective homogeneous space. 
\item[$\circ$]  A non-classical flag domain contains a compact
rational curve with ample normal bundle.

\end{itemize}
The first property  is by definition while the second
 is one of the main results of \cite{grt}.
They prove that a non-classical flag domain 
is rationally chain connected; that is,
any two points are  connected by a chain of compact rational curves contained
 in it.  The existence of an irreducible
 rational curve with ample normal bundle
follows from this by a standard smoothing argument \cite[II.7.6.1]{rc-book}.

As a simple example, $\SU(n,1)\subset \GL(n+1)$ acts on
$\p^n$ with two open orbits. One of them is the open unit ball in $\c^n$;
an Hermitian symmetric domain. The other is the complement of the closed 
unit ball; it is a non-classical flag domain. We see right away that
it contains many lines and in fact there is a conic through any two
of its points.

The following  questions   were raised by
Griffiths and  Toledo.

\begin{ques} \label{grt.ques}
 Let $X$ be a  projective, homogeneous variety
and $D\subset X$ a non-classical flag domain.
Let $L_X$ be an (algebraic) line bundle on $X$ and 
$L_D$  an (analytic) line bundle on $D$.
\begin{enumerate}
\item Is  $H^0(D, L_D)$  finite dimensional?
\item Is  the restriction map
$H^0(X,L_X)\to H^0(D,L_X|_D) $  an isomorphism?
\item Is ${\mathcal Mer}(D)={\mathcal Rat}(X)$?
\end{enumerate}
\end{ques}

Theorem \ref{main.thm} 
 answers these questions affirmatively.

We note that, by contrast, 
the   two properties marked by $\circ$   are not sufficient
to understand  higher cohomology groups, not even $H^1(D, \o_D)$.

 \begin{ack}
This paper grew out of my attempt to answer 
the above questions of P.~Griffiths and  D.~Toledo.
I also thank them
for further helpful discussions and comments.
Partial financial support   was provided  by  the NSF under grant number 
DMS-07-58275.
\end{ack}

\section{Chains of subvarieties}

\begin{say}[Chains of subvarieties in $X$] \label{chains.say}
Let $X=G/H$ be a quasi projective, homogeneous space,
$Z$ an irreducible variety and $u:Z\to X$ a morphism.
For now we are interested in the case when $u:Z\into X$ is a
subvariety, but in Section \ref{leff.sec} we use the general setting.

A {\it $Z$-chain} on length $r$ in (or over) $X$ 
consists of
\begin{enumerate}
\item points  $a_i, b_i\in Z$ for $i=1,\dots, r$ and
\item translations $\tau_i$ for $i=1,\dots, r$ such that
\item $\tau_i\bigl(u(b_i)\bigr)=\tau_{i+1}\bigl(u(a_{i+1})\bigr)$
 for $i=1,\dots, r-1$.
\end{enumerate}
The triple  $(a_i, b_i, \tau_i)$ is a {\it link} of the chain.
We also write it as
$$
\bigl(\tau_i\circ u: (Z, a_i, b_i)\to X\bigr).
$$
We say that the chain {\it starts} at $\tau_1\bigl(u(a_{1})\bigr)\in X$ 
and {\it ends}
at $\tau_r\bigl(u(b_r)\bigr)\in X$.

The points $a_i, b_i$  determine a connected, reducible variety
$Z(a_1, b_1, \dots, a_r, b_r)$
obtained from $r$ disjoint copies $Z_1,\dots, Z_r$ of $Z$ by identifying
$b_i\in Z_i$ with $a_{i+1}\in Z_{i+1}$ for $i=1,\dots, r-1$.
The morphisms $\tau_i\circ u$ then define a morphism
$$
(\tau_1\circ u,\dots, \tau_r\circ u): Z(a_1, b_1, \dots, a_r, b_r)\to X.
$$
Its image is a connected subvariety of $X$ that contains
the starting and end points of the chain and whose
irreducible components are translates of $u(Z)\subset X$. 
(For most purposes one can identify a chain with its image in $X$, but
this would be  slightly inconvenient when considering deformations of a
trivial chain where $\tau_1=\cdots=\tau_r$. The difference
becomes crucial only when we consider properties
(\ref{main.thm}.10--12).)

The set of all chains of length $r$
is naturally an algebraic  subvariety of $Z^{2r}\times G^r$.
It is denoted by $\chain(Z,r)$.
We write $\chain(Z,r,x)\subset \chain(Z,r)$ to denote the subvariety of all
chains starting at $x\in X$.
Up to isomorphism $\chain(Z,r,x)$ is independent of $x$.

The starting point (resp.\ the end point)
gives a morphism
$$
\alpha, \beta: \chain(Z,r)\to X. 
$$
Thus $\beta(\chain(Z,r,x)\bigr)\subset X$
is the set of points that can be connected to $x$ by a
$Z$-chain of length $\leq r$.

Note that $\beta(\chain(Z,r,x)\bigr)\subset X$
is constructible; let $W_r(x)$ denote its closure.
If there is a translate $\tau\bigl(u(Z)\bigr)$ that is not contained in $W_r(x)$
but whose intersection with
$W_r(x)$ is nonempty, then, by translating $\tau\bigl(u(Z)\bigr)$
 to nearby points
we see that $\dim W_{r+1}(x)>\dim W_r(x)$; see  \cite[4.13]{rc-book}.

Thus the sequence $W_1(x)\subset W_2(x)\subset \cdots$ stabilizes after at most
$\dim X$ steps with an irreducible subvariety $W(x)$.
 Furthermore, if $x'\in W(x)$ then $W(x')\subset W(x)$ hence in fact
$W(x')= W(x)$. Since $x$ and $x'$ can both be  connected
by a $Z$-chain of length $\leq \dim W(x)$ to points in a dense open subset of 
 $W(x')= W(x)$, we see that $x$ and $x'$ are connected
to each other by a $Z$-chain of length $\leq 2 \dim W(x)$.

Note also that if a $Z$-chain connects $x$ to $\tau_1(x)$ and
another one connects $x$ to $\tau_2(x)$
 then translating the second chain
and concatenating gives a $Z$-chain that 
connects $x$ to $\tau_1\bigl(\tau_2(x)\bigr)$. 

We can summarize these considerations as follows.
\end{say}

\begin{prop} \label{Z.chain.prop}
Let $X=G/H$ be a quasi projective, homogeneous space,
$Z$ an irreducible variety and $u:Z\to X$ a morphism.

Then there is a subgroup $H\subset K\subset G$
such that  two points $x_1, x_2\in X$ are connected
by a $Z$-chain iff they are contained in  the same fiber of
the natural projection $p_K:G/H\to G/K$.

Furthermore, in this case $x_1, x_2\in X$ can be connected
by a $Z$-chain of length $\leq 2(\dim K-\dim H)$. \qed
\end{prop}

\begin{say}[Equivalence of  (\ref{main.thm}.13--16)]
Let $X$ be any homogeneous space under a group $G$
and $Z\subset X$ a compact,
 irreducible, nondegenerate  subvariety. 
Thus the morphism $p_K:G/H\to G/K$ above is constant
and Proposition \ref{Z.chain.prop} implies both (\ref{main.thm}.14--15).

In order to see (\ref{main.thm}.13) let  $B\subset X$ be a nonzero divisor.
By (\ref{main.thm}.14) a suitable translate of $Z$ intersects $B$ but is not
contained in it, so $\tau_g^*\o_X(B)|_Z$ has a
nonconstant section. In particular,  $\tau_g^*\o_X(B)|_Z$
is not in $\pic^{\circ}(Z)$. Thus
$\o_X(B)|_Z$ is also not in $\pic^{\circ}(Z)$, hence it is not trivial.
Therefore
$Z\cap B\neq \emptyset$. \qed

\end{say}

\begin{say}[Ampleness of the normal bundle]
Let $X=G/H$ be a quasi projective homogeneous space
and $Z\subset X$ a degenerate subvariety.
Let $W\subset X$ be the fiber of  $p_K:G/H\to G/K$ 
as in Proposition \ref{Z.chain.prop} that contains $Z$. 
Then $I_W/I^2_W$ is a trivial bundle of rank $=\codim_XW$, hence
$$
\oplus_i \o_Z\cong \bigl(I_W/I^2_W\bigr)|_Z\into 
I_Z/I_Z^2
$$
is a trivial subsheaf of rank $=\codim_XW$.
In particular,  the normal sheaf of $Z\subset X$ is not ample in any sense.

Thus if $Z\subset X$ is a smooth (or local complete intersection) 
subvariety with ample normal bundle then $Z$ is nondegenerate.

The converse does not hold. For instance, a line in a quadric hypersurface
of dimension $\geq 3$ is nondegenerate but its normal bundle 
has a trivial summand. More generally, if $X$ is a projective homogeneous space
with Picard number 1 then a line (that is, a minimal degree rational curve)
 in $X$ has ample normal bundle 
iff $X=\p^n$.
\end{say}

\section{Proof of the finiteness conditions}

\begin{say}[Chains of subvarieties in $D$] \label{chains.D.say}
Using the notation of Paragraph \ref{chains.say}, 
note that forgetting the last component of a chain
gives a natural morphism
$\chain(Z,r+1,x)\to \chain(Z,r,x)$
whose fibers are isomorphic to $\chain(Z,1,x)$.

Furthermore,   $\chain(Z,1,x)\subset Z^2\times G$ and the
fibers of the  projection  to $Z^2$ are translates of $H$.
Thus if $H$ is connected 
then $\chain(Z,1,x)$ is irreducible
and so are the other varieties $\chain(Z,r,x)$. 

Let $\chain^0(Z,n,x)\subset \chain(Z,r,x)$
denote those chains for which $a_i, b_i\in Z$ 
(as in (\ref{chains.say}.1)) are smooth points.
This is a Zariski open condition.

For an open subset $D\subset X$, let $\chain^0(Z,D,r,x)\subset \chain^0(Z,r,x)$
 denote those chains whose image is contained in $D$.
Clearly $\chain^0(Z,D,r,x)\subset \chain(Z,r,x)$ is open
and nonempty if $x\in Z\subset D$ is a smooth point
since it contains the constant chain
where $a_i=b_i=x$ and $\tau_i=1$ for every $i$.

If $Z\subset X$ is nondegenerate then
 $\beta:\chain^0(Z,r,x)\to X$ is dominant for $r\geq 2\dim X$
and a dominant morphism is
 generically smooth. 
Thus $\beta$ is also smooth at some point of
$\chain^0(Z,D,r,x)$. We have thus established the following.
\end{say}

\begin{lem}\label{cainms.in.D.lem}
Let $X=G/H$ be a quasi projective,  homogeneous space
with connected stabilizer $H$
and $Z\subset X$  a compact, irreducible, nondegenerate subvariety.
Let $Z\subset D$ be an open neighborhood.

Then $\beta\bigl(\chain^0(Z,D,r,x)\bigr)$ contains a nonempty Euclidean open
subset  $U_r\subset  D$ for $r\geq 2\dim X$.\qed
\end{lem}

The following example shows that connectedness of $H$ is
quite important here.

\begin{exmp}\label{H.disconn.exmp}
 Let $X=S^2\p^2\setminus (\mbox{diagonal})$
with the diagonal $\GL_3$-action. Its universal cover is 
$\tilde X=\p^2\times \p^2\setminus (\mbox{diagonal})$.

Let $C_1\subset \tilde X$ be a line contained in 
some $\p^2\times \{\mbox{point}\}$ and $C\subset X$ its image.
The preimage of $C$ in $\tilde X$ is a disjoint union of a
horizontal and of a vertical line. Thus $C$-chains (of length 4) 
connect any 2 points $X$, yet $C$ has an open neighborhood
of the form $D\cong \p^1\times (\mbox{unit disc})$.
Chains of compact curves in $D$ do not connect
two general points.

In general, let $X$ be a homogeneous space and $\pi:\tilde X\to X$
its universal cover. If one (equivalently every)
irreducible component of $\pi^{-1}(Z)$ is nondegenerate,
then $Z$ has the good properties one expects based on the
simply connected case, but not otherwise.
\end{exmp}

We use  $Z$-chains in $D$ to prove that
(\ref{main.thm}.16) $\Rightarrow$ (\ref{main.thm}.1--4). 
The following lemma, modeled on \cite[Thm.2]{Nadel91},
 shows that a section that vanishes to high enough order
at one point of a $Z$-chain will vanish at all points.
If $Z$ is smooth, then one needs the semipositivity of the
normal bundle  $N_{Z,X}$; equivalently, the
seminegativity of $I_Z/I_Z^2$ where
$I_Z\subset \o_X$ is the ideal sheaf of $Z$.  If $Z$ is singular,
the seminegativity of $I_Z/I_Z^2$ alone does not seem to be enough,
one needs control of the quotients  $I_Z^{(m)}/I_Z^{(m+1)}$
of the symbolic powers; see 
Paragraphs \ref{symb.powers.defn}--\ref{pos.nbhds.lem}
 for details.

\begin{lem} \label{nadel.lem} Let $D$ be a normal complex space,
 $Z_1,\dots, Z_n\subset D$ compact subvarieties, $L$ a line bundle on
$D$ and $s\in H^0(D,L)$ a section. Assume the following.
\begin{enumerate}
\item For $j=1,\dots n$ there are smooth points $p_{j}\in  Z_{j}$
such that $p_{j}\in Z_{j-1}$ for $j\geq 2$.
\item For $j=1,\dots n$ there is a family of irreducible curves
$\{C_j(\lambda)\}$  passing through $p_j$ and covering
a dense subset of  $Z_j$ such that $\deg_{C_j(\lambda)}L\leq d_j$ for 
some $d_j\in \n$.
\item  $I_{Z_j}^{(i)}/I_{Z_j}^{(i+1)}$ are subsheaves of
a trivial sheaf $\oplus_m \o_{Z_j}$ for $i\geq 1$.
\item $s$ vanishes at $p_1$ to order  $c+\sum_{j=1}^n d_j$.
\end{enumerate}
Then $s$ vanishes along $Z_r$ to order  $c+\sum_{j=r+1}^n d_j$ for every $r$.
\end{lem}

Proof. We start with the case $i=1$ and write $Z:=Z_1$.
 Choose $q$ such that 
$$
s\in H^0\bigl(D, I_Z^{(q)}\otimes L\bigr)\setminus
H^0\bigl(D, I_Z^{(q+1)}\otimes L\bigr).
$$
Thus we get a nonzero section
$$
\bar s\in H^0\bigl(Z, \bigl(I_Z^{(q)}/I_Z^{(q+1)}\bigr)\otimes L\bigr),
$$
which vanishes at $p=p_1$ to order $c+d_1-q$.
Using assumption (3), we get at least 1 
nonzero section 
$$
\tilde s\in H^0\bigl(Z, \o_Z\otimes L\bigr)
$$
that vanishes at $p=p_1$ to order $c+d_1-q$.
Restricting this to the curves $C(\lambda)$ we see
that $\tilde s$ is identically zero on $Z$,
unless $q\geq c$.  

Returning to the general case, we see that
if $s$ vanishes at $p_1$ to order  $c+\sum_{j=1}^n d_j$
then it vanishes along $Z_1$ to order  $c+\sum_{j=2}^n d_j$,
in particular, $s$  vanishes at $p_2$ to order  $c+\sum_{j=2}^n d_j$.
Repeating the argument for the shorter chain $Z_2+\cdots+Z_n$
completes the proof. \qed

\begin{say}[Proof of (\ref{main.thm}.16) $\Rightarrow$ (\ref{main.thm}.1--4)]
Let $X=G/H$ be a quasi projective,  homogeneous space
with connected stabilizer $H$,
 $Z\subset X$  a compact, irreducible, nondegenerate subvariety and
 $D\supset Z$  an open neighborhood.
Pick a smooth point $x\in Z$.

Let $L$ be a line bundle on $D$ and  $H$  a very ample
line bundle  on $Z$. Then $Z$ is covered by a family of 
irreducible curves $\{C(\lambda)\}$ passing through $x$
obtained as intersections of $\dim Z-1$  members of $|H|$.
Set $d(L):=(L|_Z\cdot H^{\dim Z-1}\bigr)$.

We check in Lemma \ref{pos.nbhds.lem}
that  $Z$ satisfies the crucial condition (\ref{nadel.lem}.3).

By Lemma \ref{cainms.in.D.lem},  for $r\geq 2\dim X$ there is
an open subset $U_r\subset  D$ whose points can be
connected to $x$ by a $Z$-chain of length $r$ satisfying
the assumptions (\ref{nadel.lem}.1--3).

Thus if a section $s\in H^0(D,L)$  vanishes at $x$ to order 
$1+rd(L)$ then it vanishes at every point of $U_r$.
Since $U_r$ is open, this implies that
$s$  is identically zero. This shows that
$$
\dim H^0(D,L)\leq \binom{\dim X+ 2\dim(X)d(L)}{\dim X},
$$
which proves (\ref{main.thm}.2). Since $d(L^m)=md(L)$,
we also have (\ref{main.thm}.3). 

Since $H^0(D, \o_D)$ is a $\c$-algebra
without zero divisors, $H^0(D, \o_D)=\c$ is equivalent to
$\dim H^0(D, \o_D)<\infty$.

Finally consider (\ref{main.thm}.4).
We use induction on $\rank E$. If $\rank E=1$ then
its double dual $E^{**}$ is a line bundle and
$H^0(D,E)\subset H^0\bigl(D,E^{**}\bigr)$ shows that
$H^0(D,E)$ is finite dimensional. In the higher rank case, we are done if
 $H^0(D,E)=0$. Otherwise there is 
a nontrivial map $\o_D\to E$ and thus a rank 1
subsheaf $E_1\subset E$ such that $E/E_1$ is again torsion free.
Thus
$h^0(D,E)\leq  h^0(D,E_1)+h^0(D,E/E_1)$ and we are done by induction.\qed
\end{say}

\begin{defn}[Symbolic powers] \label{symb.powers.defn}
Let $X$ be a  variety and $Z\subset X$ an irreducible, reduced
subvariety with ideal sheaf $I_Z$. Let $T_m\subset \o_X/I_Z^m$
denote the largest subsheaf whose sections are supported on a 
smaller dimensional subset of
$Z$. Let $I_Z^{(m)}\subset \o_X$ denote the preimage of $T_m$.
It is called the $m$-th {\it symbolic power} of $I_Z$. 

If $X$ is smooth and  $Z$ is also smooth (or a local complete intersection)
then $I_Z^{(m)}=I_Z^m$. 

The main advantage of symbolic powers is that
the quotients  $I_Z^{(m)}/I_Z^{(m+1)}$ are torsion free sheaves on $Z$.
There are also obvious maps
$I_Z^m/I_Z^{m+1}\to I_Z^{(m)}/I_Z^{(m+1)}$ that are
isomorphisms on a dense open subset.
\end{defn}

\begin{lem} \label{pos.nbhds.lem}
Let $X$ be a homogeneous space and
$Z\subset X$ a reduced subscheme. Then
$I_Z^{(m)}/I_Z^{(m+1)}$ can be written as a subsheaf
of $\oplus_i \o_Z$. 
\end{lem}

Proof. Let us start with the $m=1$ case.
This is well known but going through it will
show the path to the general case.
Every tangent vector field $v\in H^0(X, T_X)\supset \operatorname{Lie}(G)$
gives a differentiation
$d_v:\o_X\to \o_X$ which is not $\o_X$-linear. However, if 
$\phi\in \o_X$ and $s\in I_Z$ are local sections then
$d_v(\phi\cdot  s)=d_v(\phi)\cdot  s + \phi\cdot d_v(s)$
shows that differentiation composed with restriction to $Z$
gives an $\o_X$-linear map
$d_v: I_Z\to \o_Z$. Applying this for a basis of  $ H^0(X, T_X)$
gives $I_Z/I_Z^{2}\to \oplus_i \o_Z$ whose kernel is supported
at $\sing Z$. By definition, $I_Z/I_Z^{(2)}$ has no
sections supported on a nowhere dense subset.
Thus we get an injection 
$I_Z/I_Z^{(2)}\into \oplus_i \o_Z$.

If $Z$ is smooth then $I_Z^{(m)}=I_Z^{m}$, thus
we get the required 
$$
I_Z^{(m)}/I_Z^{(m+1)}=S^m\bigl(I_Z/I_Z^{2}\bigr)
\into S^m\bigl(\oplus_i \o_Z\bigr).
$$
However, in general  
we only have a map
$$
S^m\bigl(I_Z/I_Z^{2}\bigr)\to I_Z^{(m)}/I_Z^{(m+1)}
$$
which is an isomorphism over the smooth locus. Thus
$I_Z^{(m)}/I_Z^{(m+1)} $
is more positive than $S^m\bigl(I_Z/I_Z^{2}\bigr) $.

For $I_Z^{(m)}$ we work with $m$-th order differential operators 
$D=d_{v_1}\cdots d_{v_m}$.
These give well defined maps of sheaves  $D: I_Z^{(m)}\to \o_Z$.
$\o_X$-linearity can be checked over the open set $X\setminus \sing Z$.

To simplify notation, set $M:=\{1,\dots, m\}$ and for 
 $J=\{j_1<\cdots<j_r\}\subset M$ write
$D_J:=d_{v_{j_1}}\cdots d_{v_{j_r}}$. Leibnitz-rule then says that
$$
D(\phi\cdot  s)=\tsum_{J\subset M}\
D_J(\phi)\cdot D_{M\setminus J}(s).
$$
If  $|M\setminus J|<m$  and $s\in I_Z^{(m)}$
then  $D_{M\setminus J}(s)|_Z$ vanishes on the  open set $Z\setminus \sing Z$,
hence everywhere. 
The only term left is
$D_{\emptyset}(\phi)\cdot D_{M}(s)=\phi\cdot D_{M}(s)$. 
Thus we get an $\o_X$-linear map
$$
d_{v_1}\cdots d_{v_m}: I_Z^{(m)}\to \o_Z.
$$
By letting ${v_1},\dots, {v_m}$ run through a  basis of  $ H^0(X, T_X)^{m}$,
we get the required injection
$I_Z^{(m)}/I_Z^{(m+1)}\into \oplus_i \o_Z$. \qed

\section{Meromorphic and holomorphic sections}

Here we show that property
(\ref{main.thm}.7) and  (\ref{main.thm}.13) imply
(\ref{main.thm}.5--6) in   general.

\begin{prop} Let $X$ be a normal, quasi projective variety, 
and $D\subset X$  an open subset
such that ${\mathcal Mer}(D)={\mathcal Rat}(X)$.
The following are equivalent.
\begin{enumerate}
\item $X\setminus D$ does not contain any nonzero, effective divisors.
\item Let $L$ be an ample line bundle on $X$. Then, for
every $m\in \n$,  the restriction map
$H^0(X,L^m)\to H^0(D,L^m|_D) $ is an isomorphism.
\item For every 
 reflexive, coherent sheaf $F$ on $X$,  the restriction map
$H^0(X,F)\to H^0(D,F|_D) $ is an isomorphism. 
\end{enumerate}
\end{prop}

Proof. Assume (1) and 
let $L$ be an ample line bundle on $X$
with at least one global section $s_X\neq 0$. 
Let $s_D$ be a global section of $L|_D$.
Then $s_D/s_X$ is a meromorphic function on $D$, hence, by
assumption, it extends to a rational function  $r_X$ on $X$.
Thus $r_Xs_X$ is  a rational section
 of $L$ such that $(r_Xs_X)|_D=s_D$. Since $s_D$ is holomorphic,
the polar set of $r_Xs_X$ must be disjoint from $D$.
However, $D$ meets every divisor,
 so $r_Xs_X$ has to be
a regular section of $L$.

Next we show (2) $\Rightarrow$ (3).
 Let $L$ be an ample line bundle on $X$.
 Then $F^*\otimes L^m$ is generated
by global sections for $m\gg 1$. Thus we have an injection
$j:F\into \oplus_i L^m$ of $F$ into a direct sum of many copies of $L^m$.

Since $H^0(D,F|_D)\subset\oplus_i H^0(D,L^m|_D)$,
every global section $s_D$ of $F|_D $
is the restriction of a global section  $s_X$ of $\oplus_i L^m$. 
We have two subsheaves
$$
F\subset \langle F, s_X\rangle \subset \oplus_i L^m
$$
and they agree on $D$. Thus the support
of the quotient $\langle F, s_X\rangle /F$ is disjoint from $D$. 
Since $D$ meets every divisor,
the support of  $\langle F_X, s_X\rangle /F_X$
has codimension $\geq 2$. Since $F$ is reflexive, this
forces $\langle F_X, s_X\rangle =F$, hence
$s_X\in H^0(X,F)$. 

The  converse (3) $\Rightarrow$ (2) is clear.

Finally we show  (2) $\Rightarrow$ (1). Assuming the contrary,
there is an effective  divisor $B\subset X$ that is disjoint from $D$.
Choose $m$ such that  $L^m(B)$ is generated by global sections.  
Then
$$
H^0\bigl(X, L^m\bigr)\subsetneq H^0\bigl(X, L^m(B)\bigr)\into 
H^0\bigl(D, L^m(B)|_D\bigr)=
 H^0\bigl(D, L^m|_D\bigr) 
$$
contradicts (2). \qed
\medskip

\section{Lefschetz property and meromorphic functions}

Here we show that the Lefschetz--type property
(\ref{main.thm}.9) and (\ref{main.thm}.3)
imply
(\ref{main.thm}.7) in   general.
Although we do not use it, it is worth noting that
\cite{nap-ram} proves that (\ref{main.thm}.2) implies that 
the map in (\ref{main.thm}.9) has finite index image for every $X^0$.

First we show, using (\ref{main.thm}.3) that
${\mathcal Mer}(D)$ is an algebraic extension of ${\mathcal Rat}(X)$.
Then we establish that having a
meromorphic function  on $D$ that is algebraic over ${\mathcal Rat}(X)$
is equivalent to a failure of (\ref{main.thm}.9).

\begin{prop} \cite{MR0241433, MR0232780} \label{merom.thm} 
Let $X$ be a normal, quasi projective variety of dimension $n$ and
 $D\subset X$   a Euclidean open subset.
Assume that $h^0(D, L^m)=O\bigl(m^n\bigr)$ for every line bundle $L$ on $D$.
Then ${\mathcal Mer}(D)$ is an algebraic extension of ${\mathcal Rat}(X)$.
\end{prop}

Proof. Let $f_1,\dots, f_n$ be algebraically independent
rational functions on $X$ and $\phi$  a meromorphic function on $D$.
Let $B_1,\dots, B_n$ and $B_0$ be their divisors  of poles.
Consider the line bundle
$$
L:=\o_D\bigl(B_1|_D+\cdots+ B_n|_D+B_0\bigr).
$$
We can view $f_1|_D, \dots, f_n|_D$ and $\phi$ as 
sections of $L$. Thus the monomials
$$
\bigl\{\phi^{a_0}\cdot \textstyle{\prod}_i \bigl(f_i|_D\bigr)^{a_i} :
\tsum_{i=0}^n a_i=m\bigr\}
$$
are all sections of $L^m$. 
The number of these monomials grows like  $m^{n+1}$ while, by
assumption,   the dimension of $H^0(D,L^m_D) $
grows like  $m^{n}$. Thus, for $m\gg 1$, the function
 $\phi$ satisfies a nontrivial identity
$$
\tsum_{i=0}^m\  h_i\cdot \phi^i=0\qtq{where} h_i\in  {\mathcal Rat}(X).\qed
%\eqno{(\ref{merom.thm}.2)}
$$

\begin{prop}\label{3.degs.same.prop} 
Let $X$ be a normal, quasi projective variety and
 $D\subset X$  a Euclidean open subset.
For every $d\in \n$ the following are equivalent.
\begin{enumerate}
\item There is a  $\phi\in {\mathcal Mer}(D)$  such that
$\deg \bigl[{\mathcal Rat}(X)(\phi):{\mathcal Rat}(X)\bigr]=d$.
\item There is an irreducible (possibly ramified)
cover $\pi:\tilde X\to X$ of degree $d$
such that the injection $j:D\into X$ lifts to an
injection $\tilde j:D\into \tilde X$.
\item There is a Zariski open subset $X^0\subset X$ such that
$$
\im\bigl[\pi_1\bigl(D\cap X^0\bigr)\to \pi_1\bigl(X^0\bigr)\bigr]
\subset \pi_1\bigl(X^0\bigr)
\qtq{has index $d$.}
$$
\end{enumerate}
\end{prop}

Proof. Let $\phi$ be a meromorphic function on $D$ that has degree
$d$ over $ {\mathcal Rat}(X)$. 
Let
$$
\tsum_{i=0}^d\ h_i\cdot \phi^i =0,
\eqno{(\ref{3.degs.same.prop}.4)}
$$
be the minimal polynomial of $\phi$ 
where the $h_i$ are rational functions on $X$.

Let  $\pi:\tilde X\to X$ be the
normalization of $X$ in the field ${\mathcal Rat}(X)(\phi) $.

The key observation is that  we can  
think of  $\phi$ in two new ways:
either as a rational function $\tilde \phi$ on $\tilde X$
or as 
 a multi-valued algebraic  function $\phi_X$ on $X$
whose  restriction to $D$ contains a
 single-valued branch that agrees with $\phi$.

Since (\ref{3.degs.same.prop}.4) is irreducible over ${\mathcal Rat}(X)$,
 its discriminant is not identically zero, thus
there is a dense, Zariski open subset $X^0\subset X$
such that $\pi$ is \'etale over $X^0$ and
$\tilde \phi$ 
 takes up different values at 
different points of $\pi^{-1}(x)$ 
for all  $x\in X^0$. Thus the single valued branch $\phi$ of  $\phi_X$
determines a lifting of the injection 
$$
j^0:D\cap X^0\into X \qtq{to }\tilde j^0: D\cap X^0\into \tilde X.
$$
This shows (3).

Since $\pi$ is finite, in  suitable local coordinates
we can view $\tilde j^0 $ as a bounded holomorphic function.
Thus  $\tilde j^0 $ extends to a
lifting $\tilde j: D\into \tilde X$, hence
(3) $\Rightarrow$ (2).

Conversely, it is clear that (2) $\Rightarrow$ (3).
Furthermore, if (2) holds and $\psi$ is a rational function
on $\tilde X$ that generates ${\mathcal Rat}(\tilde X)/{\mathcal Rat}(X)$
then $\psi\circ \tilde j$ is a 
 meromorphic function  on $D$ such that
${\mathcal Rat}(X)(\psi\circ \tilde j)={\mathcal Rat}(\tilde X)$, hence
$\deg \bigl[{\mathcal Rat}(X)(\psi\circ \tilde j):{\mathcal Rat}(X)\bigr]=d$.
\qed

\begin{say}[Proof of (\ref{main.thm}.7) and (\ref{main.thm}.10) $\Rightarrow$ 
(\ref{main.thm}.8)]
\label{covers.of.D.say}
With the notation and assumptions of Proposition 
\ref{merom.thm}, let $v:\bar D\to D$ be a
finite, possibly ramified, cover.
Since  ${\mathcal Mer}(\bar D)$ is an
algebraic extension of ${\mathcal Mer}( D)$, (\ref{main.thm}.7) implies that
 ${\mathcal Mer}(\bar D)$ is an algebraic extension of ${\mathcal Rat}(X)$.

The assumption (\ref{main.thm}.10)
says that 
$\pi_1\bigl(D\cap X^0\bigr)\to \pi_1\bigl(X^0\bigr)$
is surjective for every $X^0$, thus
${\mathcal Mer}(\bar D)={\mathcal Rat}(X)$
by  Proposition \ref{3.degs.same.prop}.\qed
\end{say}

\section{Lefschetz--type properties}
\label{leff.sec}

It is clear that (\ref{main.thm}.10) $\Rightarrow$  (\ref{main.thm}.9)
and (\ref{main.thm}.12) $\Rightarrow$  (\ref{main.thm}.11).

Let $u_D:\tilde D\to D$ be a finite, \'etale cover.
Those $g\in G$ for which $\tau_g(Z)\subset D$
form a Euclidean open subset of $G$, thus there are
general translations  (in the sense of (\ref{main.thm}.11--12))
 such that $\tau_g(Z)\subset D$.
There is a finite (\'etale) cover $u_z:\tilde Z\to Z$
such that $\tau_g\circ u_Z$ factors through $u_D$. 
This shows that (\ref{main.thm}.12) $\Rightarrow$ (\ref{main.thm}.10).

It remains to show that if $Z$ is nondegenerate then (\ref{main.thm}.12) holds.
Thus let $u_Z:\tilde Z\to Z$
be a finite  cover  and $X^0\subset X$ a
Zariski open subset.
For this  we use $\tilde Z$-chains  in
 general position.

\begin{say}[Chains  in general position]\label{gen.pos.trans.say}
Let $m: G\times \tilde Z\to X$ be the $G$-action composed with $u_Z$.
Every map between algebraic varieties is a 
locally topologically trivial fiber bundle
over a Zariski open subset, cf.\ \cite[p.43]{gm-book}. Thus
there is a Zariski open subset  $G^0\subset G$ such that
the first projection
$\pi_G: G\times \tilde Z\to G$
restricts to a topologically trivial fiber bundle
on $m^{-1}(X^0)$.
We denote its fibers by
$\tilde Z^0_g:=(\tau_g\circ u_Z)^{-1} (X^0)$.
Marking  a pair of  smooth points has no significant  effect topologically, 
thus the triples
$\bigl(a, b, \tilde Z^0_g\bigr)$, where  $a\neq b$ 
are smooth points of $\tilde Z^0_g $,
are fibers of a topologically trivial fiber bundle
over a  Zariski open subset  of $G^0\times \tilde Z^2$.

We say that $(a, b ,\tau_g)\in \chain(\tilde Z,1,x)$
is in {\it general position}  with respect to $X^0$
if $g\in G^0$, $a\neq b$ are smooth points of $Z$ and they are both mapped to
$X^0$. 

The  set of all general position maps forms a
Zariski open subset 
 $\chain^*(\tilde Z,1,x)\subset \chain(\tilde Z,1,x)$
which is nonempty for general $x\in X$. 

For us  a key point is that 
 the image of the induced map
 $$
\Gamma(X^0,x)
:=\im\bigl[\pi_1\bigl(\tilde Z^0_g, a\bigr)\to \pi_1\bigl( X^0, x\bigr)\bigr]
\subset \pi_1\bigl( X^0, x\bigr)
\eqno{(\ref{gen.pos.trans.say}.1)}
 $$
is independent of 
$(a, b, \tau_g)\in \chain^*(\tilde Z,1,x)$
whenever the latter is nonempty.

We say that a  $\tilde Z$-chain as in (\ref{chains.say}.1--3) is in 
{\it general position} with respect to $X^0$
if every link  $(a_i, b_i, \tau_{g_i})$ is in general position.

As before,  $\tilde Z$-chains  in 
 general position with respect to $X^0$ form a 
Zariski open subset
$\chain^*(\tilde Z,r,x)\subset \chain(\tilde Z,r,x)$
which is nonempty for general $x\in X$.

Corresponding to
$\Gamma(X^0,x)\subset  \pi_1\bigl(X^0,x\bigr)$ as in (\ref{gen.pos.trans.say}.1)
there is an \'etale cover
$$
\pi_X^0: \bigl(\tilde X^0, \tilde x\bigr)\to \bigl(X^0,x\bigr)
\eqno{(\ref{gen.pos.trans.say}.2)}
$$
such that every general position map
$\tau_g\circ u_Z: \bigl(\tilde Z^0_g,a\bigr)\to \bigl(X^0,x\bigr)$ 
lifts to 
$$
\widetilde{\tau_g\circ u_Z}: \bigl(\tilde Z^0_g,a\bigr)\to 
\bigl(\tilde X^0,\tilde x\bigr).
\eqno{(\ref{gen.pos.trans.say}.3)}
$$
(We do not yet know that $\Gamma(X^0,x)$ has finite index,
so $\tilde X^0\to X^0 $ could be an infinite degree cover.)
Note further that (\ref{gen.pos.trans.say}.1) implies the following.
\medskip

{\it Claim} \ref{gen.pos.trans.say}.4. Let $\tilde x'\in \tilde X^0$
be any preimage of a point  $x'\in X^0$. Assume that  
$$
\tau'\circ u_Z: \bigl(\tilde Z^0_g,a'\bigr)\to \bigl(X^0,x'\bigr)
\qtq{lifts to} 
\widetilde{\tau'\circ u_Z}: \bigl(\tilde Z^0_g,a'\bigr)\to 
\bigl(\tilde X^0,\tilde x'\bigr)
$$
for some $(a', b', \tau')\in \chain^*(\tilde Z,1,x')$. Then the lift exists
for every $(a, b, \tau)\in \chain^*(\tilde Z,1,x')$.\qed
\medskip

\end{say}

\begin{prop}\label{key.lift.prop.1}
Every $\tilde Z$-chain in general position with respect to $X^0$
and starting at $x$ lifts to a 
$\tilde Z$-chain on $\tilde X^0 $  starting at $\tilde x$.
\end{prop}

Proof. A  $\tilde Z$-chain is given by the data
$(a_i, b_i, \tau_i)$. By the choice of $\Gamma(X^0,x)$,
$$
\tau_1\circ u_Z: (\tilde Z^0_1, a_1)\to (X, x_1)
\qtq{lifts to}
\widetilde{\tau_1\circ u_Z}: (\tilde Z^0_1, a_1)\to (\tilde  X^0, \tilde x_1).
$$
If we let $\tilde x_2$ denote the image of $b_1$ 
then we can view the latter map as
$$
\widetilde{\tau_1\circ u_Z}: (\tilde Z^0_1, b_1)\to (\tilde  X^0, \tilde x_2).
$$
Both
$$
\tau_1\circ u_Z: ( \tilde Z^0_1, b_1)\to (X, x_2)\qtq{and}
\tau_2\circ u_Z: ( \tilde Z^0_2, a_2)\to (X, x_2)
$$
are in general position with respect to $X^0$, thus
by (\ref{gen.pos.trans.say}.4),
if one of them lifts to 
$(\tilde  X^0, \tilde x_2)$ then so does the other.
This gives us
$$
\widetilde{\tau_2\circ u_Z}: (\tilde Z^0_2, a_2)\to (\tilde  X^0, \tilde x_2).
$$
We can iterate the argument to lift the whole chain. \qed
\medskip

For $r=2\dim X$, we thus get a lift of the end point map
$$
\beta_r: \chain^*(\tilde Z,r,x)\to X^0
\qtq{to}
\tilde \beta_r: \chain^*(\tilde Z,r,x)\to \tilde X^0.
$$
Since $\beta_r$ is dominant, the induced map
$\pi_1\bigl(\chain^*(\tilde Z,r,x)\bigr)\to \pi_1\bigl(X^0\bigr)$
has finite index image; cf.\  \cite[2.10]{shaf-book}. 
Thus $  \tilde  X^0$ is a finite degree cover of $X^0$ and so
 it uniquely  extends to a finite ramified cover
$\pi_X:\tilde X\to X$ where $\tilde X$ is normal.

If $\tilde X= X$ then 
$\Gamma(X^0,x)=\pi_1\bigl( X^0, x\bigr)$, thus
$\pi_1\bigl(\tilde Z^0_g, a\bigr)\to \pi_1\bigl( X^0, x\bigr)$
is surjective. This proves (\ref{main.thm}.12).

All that remains is to derive a contradiction
if $\tilde X\neq  X$.
Since $X$ is simply connected,  in this case $\tilde X\to X$
has a nonempty branch divisor  $B\subset X$. 
We use the branch divisor to show that some
chains do not lift, thereby arriving at a contradiction.

\begin{prop}\label{key.lift.prop.2}
If $\tilde X\neq  X$ then there is a 
 $\tilde Z$-chain in general position with respect to $X^0$
and starting at $x$ that does not lift to a 
$\tilde Z$-chain on $\tilde X^0 $  starting at $\tilde x$.
\end{prop}

Proof. Set $d:=\deg \tilde X/ X$.
 Let $\tau_{r+1}\circ u_Z: \tilde Z\to X$ be a general translate
whose image intersects the branch divisor $B$ generically transversally;
that is,  the scheme theoretic inverse image
$(\tau_{r+1}\circ u_Z)^{-1}(B)\subset \tilde Z $
is a reduced divisor.  Then the pull-back
$$
\tilde Z\times_X\tilde X\to \tilde Z
$$
is a degree $d$  cover that ramifies over at least 1 smooth
point of $\tilde Z$. The cover need not be connected or normal, but,
due to the ramification, 
it can not be a  union of $d$ trivial covers
$\tilde Z\cong \tilde Z$.

Let $a\in \tilde Z$ be a general smooth point
and $\tilde a_1,\dots, \tilde a_d$ its preimages.
Thus, for at least one $\tilde a_i$, the identity map
$( \tilde Z, a)\to (\tilde Z, a)$ can not be lifted to
$( \tilde Z, a)\to\bigl( \tilde Z\times_X\tilde X, \tilde a_i\bigr)$.
Thus if $x\in X$ is the image of $a$ and
$\tilde x_1,\dots, \tilde x_d\in \tilde X$ its preimages,
then for at least one $\tilde x_i$, the map
$\tau_{r+1}\circ u_Z$ can not be lifted to
$$
\widetilde{\tau_{r+1}\circ u_Z}:( \tilde Z, a)\not\to 
\bigl( \tilde X, \tilde a_i\bigr).
$$
Consider now the dominant map
$\tilde \beta_r: \chain^*(\tilde Z,r,x)\to \tilde X^0$
and let $X^*\subset X^0$ be a Zariski open subset
such that $\pi_X^{-1}(X^*)\subset \im \tilde \beta_r$.
By choosing the above $\tau_{r+1}\circ u_Z: \tilde Z\to X$
generally, we may assume that there is a smooth point
$a_{r+1}\in Z$ such that $x^*:=(\tau_{r+1}\circ u_Z)(a_{r+1})\in X^*$.

Thus, for every $\tilde x^*_i\in \pi_X^{-1}(x^*)$
there is a $\tilde Z$-chain of length $r$
whose lift to $\tilde X$ connects $\tilde x$ and $\tilde x^*_i$. 
We can add $\tau_{r+1}\circ u_Z: \bigl(\tilde Z, a_{r+1}, b_{r+1}\bigr)\to X$ 
as the last link of
any of these chains.  Thus we get $d$ different 
$\tilde Z$-chains of length $r+1$
and at least one of the can not be lifted to $\tilde X$.\qed
\medskip

This completes the proof of the
last implication (\ref{main.thm}.16) $\Rightarrow$  (\ref{main.thm}.12).

%\bibliography{refs-main/refs}

\def\cprime{$'$} \def\cprime{$'$} \def\cprime{$'$} \def\cprime{$'$}
  \def\cprime{$'$} \def\cprime{$'$} \def\dbar{\leavevmode\hbox to
  0pt{\hskip.2ex \accent"16\hss}d} \def\cprime{$'$} \def\cprime{$'$}
  \def\polhk#1{\setbox0=\hbox{#1}{\ooalign{\hidewidth
  \lower1.5ex\hbox{`}\hidewidth\crcr\unhbox0}}} \def\cprime{$'$}
  \def\cprime{$'$} \def\cprime{$'$} \def\cprime{$'$}
  \def\polhk#1{\setbox0=\hbox{#1}{\ooalign{\hidewidth
  \lower1.5ex\hbox{`}\hidewidth\crcr\unhbox0}}} \def\cdprime{$''$}
  \def\cprime{$'$} \def\cprime{$'$} \def\cprime{$'$} \def\cprime{$'$}
\providecommand{\bysame}{\leavevmode\hbox to3em{\hrulefill}\thinspace}
\providecommand{\MR}{\relax\ifhmode\unskip\space\fi MR }
% \MRhref is called by the amsart/book/proc definition of \MR.
\providecommand{\MRhref}[2]{%
  \href{http://www.ams.org/mathscinet-getitem?mr=#1}{#2}
}
\providecommand{\href}[2]{#2}

\vskip1cm

\noindent Princeton University, Princeton NJ 08544-1000

{\begin{verbatim}kollar@math.princeton.edu\end{verbatim}}

\end{document}